\documentclass[12pt]{article}
\def\Box{\framebox[2.45mm]{}}
\newtheorem{theorem}{Th\'eor\`eme}
\newtheorem{coro}[theorem]{Corollaire}

\newenvironment{remarque}{\noindent{\bf Remarque. }}{\smallskip}
\newenvironment{proof}{{\em D\'emonstration. }}{\mbox{}\hfill 
$\Box$\medskip}
\begin{document}
\title{G\'en\'eralisation de formules de type Waring}
\frenchspacing
\author{Fr\'ed\'eric Jouhet et Jiang Zeng\\
\small Institut Girard Desargues\\
\small Universit\'e Claude Bernard (Lyon 1)\\
\small 69622 Villeurbanne Cedex, France\\
\small Email: jouhet@desargues.univ-lyon1.fr\\
 \small \phantom{Email: } zeng@desargues.univ-lyon1.fr
}
\date{}
\maketitle
\begin{abstract}
We evaluate the symmetric functions $e_k$, $h_k$ and $p_k$ on the  
alphabet $\{x_r/(1-tx_r)\}$ by elementary methods and give the related generating functions.
Our formulas lead to  a new and short proof of an ex-conjecture of Lassalle \cite{Las2}, which was
proved by  Lascoux and Lassalle \cite{LL} in the framework of
 $\lambda$-rings theory. 
\end{abstract}

\section{Introduction}            
L'un des probl\`emes fondamentaux dans l'\'etude de 
fonctions sym\'etriques est le d\'eveloppement d'une 
fonction sym\'etrique sur certaines bases lin\'eaires de 
l'alg\`ebre des fonctions sym\'etriques.
Un r\'esultat classique de Waring  explicite le d\'eveloppe\-ment des 
fonctions sym\'e\-triques puissances $p_n$ dans la base lin\'eaire des
fonctions sym\'etriques \'el\'ementaires $(e_\lambda)$.

Dans cet article nous g\'en\'eralisons
la formule de Waring en d\'eveloppant les fonctions sym\'e\-triques
puissances $p_n$ \'evalu\'ees sur l'alphabet $Y=\left\{{x_1/(1-tx_1)}\right.$, 
$\left.{x_2/(1-tx_2)},\ldots\right\}$ dans la base
lin\'eaire des fonctions sym\'etriques \'el\'ementaires $(e_\lambda)$ \'evalu\'ees sur l'alphabet $X=\{x_1,x_2,\ldots\}$. De m\^eme nous consid\'erons
le probl\`eme inverse, c'est-\`a-dire, le  d\'eveloppement des
fonctions sym\'etriques $h_n$ et $e_n$ \'evalu\'ees sur l'alphabet $Y$
dans la base des fonctions puissances $p_\mu$ \'evalu\'ees
 sur l'alphabet $X$. Dans le dernier cas nous aurons besoin
d'un coefficient binomial g\'en\'eralis\'e introduit par Lassalle \cite{Las1}.
Nous en
d\'eduisons ensuite, comme applications, des d\'eveloppements int\'eressants, 
qui conduisent en
particulier à de nouvelles preuves des
ex-conjectures de Lassalle~\cite{Las1,Las2}.
D'autres preuves de ces conjectures ont \'et\'e tout r\'ecemment donn\'ees par
Lascoux et Lassalle~\cite{LL} dans le cadre des $\lambda$-anneaux. 
Notre approche repose essentiellement sur
l'op\'erateur diff\'erentiel de l'alg\`ebre des s\'eries formelles.
Il est  remarquable que l'\'etude d'un probl\`eme 
si \'el\'ementaire puisse conduire \`a une preuve tr\`es 
simple de l'identit\'e de Lascoux et Lassalle.

Nous terminons cette introduction par un  rappel \cite[Chap.1]{Mac}  
des formules qui seront utilis\'ees
dans la suite.
Observons d'abord que
\begin{equation}\label{eq:keylem1}
\sum_{n\geq 1}{n-1\choose k-1}a_nt^{n-1}=
{t^{k-1}\over (k-1)!}{d^{k-1}\over dt}\left(\sum_{n\geq 
1}a_nt^{n-1}\right).
\end{equation}
Comme 
les fonctions puissances $p_n(X) =\sum_{r\geq 1}x_r^n$ satisfont
$\sum_{n\geq 1}p_n(X)t^{n-1}=\sum_{r\geq 1}{x_r/(1-x_rt)}$,
et pour tout $k\geq 1$
$$
{d^{k-1}\over dt}\left({1\over 1-xt}\right)=(k-1)!{x^{k-1}\over (1-xt)^{k}},
$$
nous en d\'eduisons donc
\begin{equation}\label{eq:keylem2}
{d^{k-1}\over dt}\left(\sum_{n\geq 1}p_n(X)t^{n-1}\right)={(k-1)!\over 
t^k}
p_k\left({tx_1\over 1-tx_1}, {tx_2\over 1-tx_2} ,\ldots\right).
\end{equation}
Pour toute partition d'entiers $\mu$ on pose  
$
z_\mu=\prod_{i\geq 1}i^{m_i(\mu)}m_i(\mu)!,
$
o\`u $m_i(\mu)$ est le nombre de parts dans $\mu$ 
\'egales \`a $i\geq 1$, et pour tout entier $n$ positif on d\'efinit le coefficient 
binomial g\'en\'eralis\'e $\langle{\mu\atop n}\rangle$ comme \'etant le nombre de fa\c cons 
de choisir $n$ \'el\'ements dans le 
diagramme de Ferrers de $\lambda$, dont au moins un par ligne. 

Les fonctions sym\'etriques $h_n(X)$
et  $e_n(X)$ sont li\'ees aux fonctions puissances
$p_\mu(X)=\prod_{r\geq 1}p_{\mu_r}(X)$ par la formule :
\begin{eqnarray}\label{eq:waring1}
h_n(X)&=&\sum_{\mu \vdash n}z_\mu^{-1}p_\mu(X),\\
e_n(X)&=&\sum_{\mu \vdash n}(-1)^{n-l(\mu)}z_\mu^{-1}p_\mu(X).
\label{eq:waring2}
\end{eqnarray}
L'inverse de la derni\`ere est appel\'ee \emph{formule de Waring}
\cite{Macmahon}:
\begin{equation}\label{eq:waring1}
p_n(X)=\sum_{\lambda\vdash n}(-1)^{n-l(\lambda)}
{n(l(\lambda)-1)!\over \prod_im_i(\lambda)!}e_\lambda(X).
\end{equation}
Par l'involution $\omega$ d\'efinie par $\omega(e_n)=h_n$
on a aussi \cite[p. 24]{Mac}
\begin{equation}\label{eq:waring2}
p_n(X)=\sum_{\lambda\vdash n}(-1)^{l(\lambda)-1}
{n(l(\lambda)-1)!\over \prod_im_i(\lambda)!}h_\lambda(X).
\end{equation}
On note $m_\mu(X)$ la fonction sym\'etrique monomiale 
associ\'ee \`a la partition $\mu$. 
  
Nous remercions \textsc{Michel Lassalle} pour ses remarques amicales sur une 
version ant\'erieure de cet article.
\section{R\'esultats principaux}
Soit $X=\{x_1,x_2,\ldots\}$ un ensemble fini ou infini 
d'ind\'eter\-min\'ees et ${X\over 1-tX}$ l'alphabet
$\{{x_1\over
1-tx_1},{x_2\over 1-tx_2},\ldots\}$.
\begin{theorem} Pour tout $k\geq 1$ on a 
\begin{eqnarray}
p_k\left({X\over 1-tX}\right)&=&\sum_{|\mu|\geq k}t^{|\mu|-k}
{|\mu|\choose k}(-1)^{|\mu|-l(\mu)}\frac{k\,(l(\mu)-1)!}
{\prod_im_i(\mu)!}e_\mu(X),\label{eq:p1}\\
p_k\left({X\over 1-tX}\right)&=&\sum_{|\mu|\geq k}t^{|\mu|-k}
{|\mu|\choose k}(-1)^{l(\mu)-1}\frac{k\,(l(\mu)-1)!}
{\prod_im_i(\mu)!}h_\mu(X).\label{eq:p2} 
\end{eqnarray}
\end{theorem}
\begin{proof}
Les formules  (\ref{eq:keylem1}) et (\ref{eq:keylem2}) impliquent directement 
$$
p_k\left({X\over 1-tX}\right)=\sum_{j\geq k}t^{j-k}{j-1\choose k-1}p_j(X).
$$
On en d\'eduit donc (\ref{eq:p1}) et (\ref{eq:p2}) respectivement de
(\ref{eq:waring1}) et (\ref{eq:waring2}).
\end{proof}

Par la m\^eme m\'ethode nous obtenons le r\'esultat suivant.
\begin{theorem}  Pour tout entier $k\geq 1$ on  a
\begin{eqnarray}
h_k\left({X\over 1-tX}\right)&=&\sum_{|\mu|\geq k}t^{|\mu|-k}
\frac{\langle{\mu\atop k}\rangle}{z_\mu}p_\mu(X),\label{eq:h}\\
e_k\left({X\over 1-tX}\right)&=&\sum_{|\mu|\geq k}t^{|\mu|-k}
(-1)^{k-l(\mu)}\frac{\langle{\mu\atop k}\rangle}{z_\mu}
p_\mu(X).\label{eq:e}
\end{eqnarray}
\end{theorem}
\begin{proof} Notons d'abord que
$$
\left\langle{\mu\atop k}\right\rangle=\sum_{k_1+\cdots 
+k_l=k\atop k_1,\ldots, k_l\geq 1}\prod_{i=1}^l{\mu_i\choose k_i},
$$
o\`u $l=l(\mu)$.
Comme chaque partition $\mu$ de $j$ correspond \`a
$l(\mu)!/\prod_{i\geq 1}m_i(\mu)!$ compositions $(k_1,\ldots, k_l)$ de $j$ telles que
$(k_1, \ldots,k_l)$ soit une permutation des parts de $\mu$,  nous avons, 
en tenant compte de  (\ref{eq:keylem1})
 et (\ref{eq:keylem2}), 
\begin{eqnarray*}
&&\sum_{j\geq k}t^{j-k}\sum_{\mu\vdash j}
\frac{\alpha^{k-l(\mu)}}{z_\mu}\left\langle{\mu\atop k}\right\rangle p_\mu(X)\\
&&\hskip 1 cm =\sum_{k_1+\cdots +k_l=k\atop l\geq 1}{\alpha^{k-l}t^{-k}\over l!k_1\cdots k_l}
\prod_{r=1}^l\sum_{\mu_r\geq 1}{\mu_r-1\choose k_r-1}p_{\mu_r}(X)t^{\mu_r}\\
&&\hskip 1 cm =\sum_{k_1+\cdots +k_l=k\atop l\geq 1}{\alpha^{k-l}\over l!{k_1}!\cdots 
{k_l}!}
\prod_{r=1}^l{d^{k_r-1}\over dt}\left(\sum_{n\geq 
1}p_n(X)t^{n-1}\right)\\
&&\hskip 1 cm =\sum_{\mu\vdash k}\frac{\alpha^{k-l(\mu)}}{z_\mu}p_\mu\left({x_1\over 
1-tx_1},{x_2\over
1-tx_2},\ldots\right).
\end{eqnarray*}
En posant $\alpha=1$ (resp. $-1$),
nous en d\'eduisons (\ref{eq:h}) (resp. (\ref{eq:e})) en appliquant  (\ref{eq:waring1}) 
(resp. (\ref{eq:waring2})).
\end{proof}

\begin{remarque} 1) Lorsque $t=0$ on retrouve les formules classiques de type
Waring.\\
2) Dans les th\'eor\`emes 1 et 2, $t$ n'est qu'un param\`etre d'homog\'en\'eit\'e,
 mais vu le r\^ole important qu'il joue dans notre  
d\'emonstration, nous pr\'ef\'erons garder 
cette forme.
\end{remarque}

Rappelons que $h_n(X)$ et $e_n(X)$ ont pour fonctions g\'en\'eratrices:
\begin{eqnarray}
\sum_{n\geq 0}h_n(X)t^{n}&=&\prod_{r\geq 1}{1\over 1-x_rt},\label{eq:waring1'}\\
\sum_{n\geq 0}e_n(X)t^n&=&\prod_{i\geq 1}\left(1+x_it\right).\label{eq:waring2'}
\end{eqnarray}

\begin{theorem}\label{keylem} Soit $z$ et $X=\{x_1,x_2,\ldots\}$ des ind\'etermin\'ees ind\'ependantes.
Alors la s\'erie formelle 
$$
F(t,u)=(1+u)^z\prod_{r\geq 1}\left(1+{u\over 1+u}\, {tx_r\over 1-tx_r}\right)
$$
admet les trois d\'eveloppements suivants
\begin{eqnarray}
F(t,u)& =&\sum_{i,j\geq 0}u^it^j\sum_{l(\mu) \le i,|\mu|=j} {z-l(\mu)\choose i-l(\mu)}m_\mu (X),\label{eq:monome}\\
F(t,u)&=&\sum_{i,j\geq 0}u^it^j\sum_{k=0}^{\min(i,j)}{z-j\choose i-k}
\sum_{\mu\vdash j}\frac{\langle{\mu\atop k}\rangle}{z_\mu}p_\mu(X)\label{eq:puissance1},\\
 F(t,u)&=&\sum_{i,j\geq 0}u^it^j\sum_{k\geq 0}{z-k\choose i-k}
\sum_{\mu\vdash j}(-1)^{k-l(\mu)}\frac{\langle{\mu\atop 
k}\rangle}{z_\mu}p_\mu(X).\label{eq:puissance2}
\end{eqnarray}
\end{theorem}
\begin{proof} 
Tout d'abord, par d\'efinition nous avons
\begin{eqnarray*}
F(t,u)
&=&\sum_{k\geq 0}u^k(1+u)^{z-k}
\sum_{1\leq r_1<r_2<\cdots<r_k\atop m_1,\ldots, m_k\geq 1}(tx_{r_1})^{m_1}\cdots 
(tx_{r_k})^{m_k}\\
&=&\sum_{i,j\geq 0}u^it^j\sum_{k\geq 0} {z-k\choose i-k}\sum_{l(\mu)=k,|\mu|=j}m_\mu (X).
\end{eqnarray*}
D'o\`u (\ref{eq:monome}). Ensuite, dans le membre de droite de  (\ref{eq:puissance1})
en rempla\c cant $i$ par $i+k$, nous obtenons en appliquant la formule du bin\^ome 
$$
\sum_{k\geq 0}u^k
\sum_{j\geq 0}(1+u)^{z-j}t^j\sum_{\mu\vdash j}\frac{\langle{\mu\atop k}\rangle}{z_\mu}p_\mu(X),
$$
qui s'\'ecrit, en posant $s=t/(1+u)$ et en appliquant (\ref{eq:h}) et (\ref{eq:waring1'}),
$$
(1+u)^z\sum_{k\geq 0}u^kh_k\left({sx_1\over 1-sx_1},
{sx_2\over 1-sx_2},\ldots\right)=(1+u)^z\prod_{r\geq 1}\left(1-{usx_r\over 1-sx_r}\right)^{-1}.
$$
Ceci est clairement \'egal \`a $F(t,u)$.
Enfin nous d\'eduisons (\ref{eq:puissance2}) de fa\c con analogue en appliquant
(\ref{eq:e}) et (\ref{eq:waring2'}).
\end{proof}

\begin{coro}
Soit $z$ et $X=\{x_1,x_2\ldots\}$ des ind\'etermin\'ees ind\'ependantes.
Pour tous  entiers $i,j\geq 1$ on a
\begin{eqnarray*}
\sum_{l(\mu) \le i,|\mu|=j} {z-l(\mu)\choose i-l(\mu)}m_\mu (X)&=&
\sum_{k=0}^{\min(i,j)}
{z-j\choose i-k}\sum_{\mu\vdash j}
\frac{\langle{\mu\atop k}\rangle}{z_\mu}p_\mu(X)\\
&=&\sum_{k=0}^{\min(i,j)}
{z-k\choose i-k}\sum_{\mu\vdash j}(-1)^{k-l(\mu)}
\frac{\langle{\mu\atop k}\rangle}{z_\mu}p_\mu(X).
\end{eqnarray*}
\end{coro}
Comme $(p_{\mu})_{\mu}$ forme une base lin\'eaire de l'alg\`ebre des fonctions
sym\'etriques, on d\'eduit du corollaire~3 le r\'esultat suivant.
\begin{coro} Soit $z$ une variable. 
Pour des entiers $i,j\geq 1$ et toute partition $\mu\vdash j$ on a
$$
\sum_{k=0}^{min(i,j)}{z-j\choose i-k}\left\langle{\mu\atop k}\right\rangle
= \sum_{k=0}^{min(i,j)}(-1)^{k-l(\mu)}{z-k\choose 
i-k}\left\langle{\mu\atop
k}\right\rangle.
$$
\end{coro}
Enfin le corollaire~3 implique aussi le r\'esultat suivant, d\^u \`a Lascoux-Lassalle~ \cite[Lemme~2]{LL}.
\begin{coro} Pour tous  entiers $k,j\geq 1$ on a
$$\sum_{l(\mu)=k,|\mu|=j}m_\mu (X)=\sum_{\mu\vdash j}(-1)^{k-l(\mu)}
\frac{\langle{\mu\atop k}\rangle}{z_\mu}p_\mu(X).
$$
\end{coro}

 \begin{remarque} On trouvera d'autres  formules sur
la somme $\sum_{l(\mu)=k,|\mu|=j}m_\mu (X)$ dans
Macdonald \cite[p. 33 et 68]{Mac}.
\end{remarque}           
\section{Applications}
On identifie chaque  partition $\lambda$ avec son \emph{diagramme de
Ferrers} et on pose
$$
(x)_\lambda=\prod_{(i,j)\in \lambda}\left(x+j-1-(i-1)/\alpha\right).
$$
Lorsque $\lambda=(n)$ est une partition-ligne on retrouve la d\'efinition 
habituelle de factorielle montante $(x)_n=x(x+1)\cdots (x+n-1)$. 
 Par un calcul direct et en posant $Z=\left\{j-1-(i-1)/\alpha\right\}$ (${(i,j)\in \lambda}$), nous obtenons :
\begin{eqnarray*}
\frac{(y-x)_\lambda}{(y)_\lambda}&=& 
\prod_{z\in Z}\left(1-\frac{x/y}{1+z/y}\right)\\
&=&\sum_{i=0}^{|\lambda|}(-x/y)^i\sum_{z_1,\ldots, z_i\in Z}
\frac{1}{1+z_1/y}\cdots \frac{1}{1+z_i/y}\\
&=&\sum_{i=0}^{|\lambda|}\sum_{j=0}^\infty(-1)^{i+j}\frac{x^i}{y^{i+j}}
\sum_{l(\mu) \le i,|\mu|=j} {|\lambda|-l(\mu)\choose i-l(\mu)}m_\mu (Z).
\end{eqnarray*}
Nous d\'eduisons donc du corollaire~3 une courte preuve d'un r\'esultat de
Lascoux-Lassalle \cite[Thm.~4]{LL}, qui fut conjectur\'e par Lassalle
\cite[Conj.~2]{Las1}.
\begin{theorem} Soient $x$, $y$ deux ind\'etermin\'ees ind\'ependantes. 
Pour toute partition $\lambda$ soit $X=\left\{j-1-(i-1)/\alpha\right\}$, ${(i,j)\in \lambda}$, alors
$$
\frac{(y-x)_{\lambda}}{(y)_{\lambda}}=\sum_{i=0}^{|\lambda|} 
\sum_{j=0}^{+\infty}
(-1)^{i+j}\frac{x^i}{y^{i+j}}\sum_{k=0}^{\min(i,j)}
{|\lambda|-j\choose i-k}\sum_{\mu\vdash j}
\frac{\langle{\mu\atop k}\rangle}{z_\mu}p_\mu(X).
$$
\end{theorem}
En fait Lascoux et Lassalle \cite{LL} ont d\'eduit le th\'eor\`eme~5 d'un
r\'esultat plus g\'en\'eral, qui  fut  aussi conjectur\'e par Lassalle \cite{Las2}.
Nous en donnons aussi une nouvelle preuve.
\begin{theorem}
Soient $z$,$u$ et $X=\{x_1,x_2\ldots\}$ des ind\'etermin\'ees 
ind\'ependantes.
Pour tous entiers $n,\, r\geq 1$ on a
\begin{eqnarray*}
&&\sum_{\mu\vdash n}\frac{(-1)^{r-l(\mu)}}{z_\mu}\left\langle{\mu\atop 
r}\right\rangle
\prod_{i\geq 1}\left(z+\sum_{k\geq 1}u^k{(i)_k\over 
k!}x_k\right)^{m_i(\mu)}=\\
&&\hskip 2cm \sum_{j\geq 0}u^j{n+j-1\choose n-r}\sum_{k=0}^{\min(r,j)}
{z-j\choose r-k}\sum_{\mu\vdash j}
\frac{\langle{\mu\atop k}\rangle}{z_\mu}\prod_{i\geq 1}x_i^{m_i(\mu)}.
\end{eqnarray*}
\end{theorem}
\begin{proof}
Soit  $Y=\{y_1,y_2,\ldots \}$ une famille infinie d'ind\'etermin\'ees.
Comme les fonctions puissances $p_i(Y)$ sont alg\'ebriquement 
ind\'ependantes dans ce cas,
 nous pouvons supposer $x_i=p_i(Y)$ pour $i\geq 1$. 
En multipliant le membre de gauche par $t^nq^r$ et sommant sur $n,\, r\geq 1$
nous pouvons \'ecrire  sa fonction g\'en\'eratrice comme suit (voir l'Appendice ci-apr\`es):
\begin{equation}\label{eq:las}
F(tu, Tq)=(1+Tq)^{z}\prod_{j\geq 1}
\left(1+{Tq\over 1+Tq}\, {Tuz_j\over 1-Tuz_j}\right),
\end{equation}
o\`u $T=t/(1-t)$ et $ z_j={y_j/t}$ pour $j\geq 1$.
Nous en  d\'eduisons  par l'application du  th\'eor\`eme \ref{keylem} que
\begin{eqnarray*}
F(tu, Tq)&=&\sum_{r,j,k\geq 0}T^{r+j}q^{r}u^j{z-j\choose r-k}
\sum_{\mu\vdash j}\frac{\langle{\mu\atop k}\rangle}{z_\mu}p_\mu(Z)\\
&=&\sum_{r,j,k\geq 0}{t^{r}q^{r}u^j\over (1-t)^{r+j}}
{z-j\choose r-k}\sum_{\mu\vdash j}
\frac{\langle{\mu\atop k}\rangle}{z_\mu}p_\mu(Y).
\end{eqnarray*}
En \'ecrivant
$$
{t^{r}\over (1-t)^{r+j}}=\sum_{n\geq r}{n+j-1\choose n-r}t^n,
$$
nous remarquons que l'expression plus haut est aussi la fonction g\'en\'eratrice du membre de droite. 
\end{proof}

\section*{Appendice. Calcul de la fonction g\'en\'eratrice}
Afin de rendre la lecture autonome nous incluons ici une preuve classique de (\ref{eq:las}). 
Remarquons d'abord que pour toute partition $\mu$
$$
\sum_{r\geq 1}\left\langle{\mu\atop r}\right\rangle q^r=\prod_{i\geq
1}\left((1+q)^i-1\right)^{m_i(\mu)},
$$
et que  la formule du bin\^ome $
(1-x)^{-\alpha}=\sum_{n\geq 0}x^n{(\alpha)_n/n!}$ permet 
d'\'ecrire
$$
\sum_{n\geq 1}u^n\frac{(i)_n}{n!}p_n(Y)=\sum_{j\geq 1}\sum_{n\geq 
1}u^n\frac{(i)_n}{n!}y_j^n=
\sum_{j\geq 1}((1-y_ju)^{-i}-1).
$$
En multipliant le membre de gauche par $t^nq^r$ et sommant sur $n,\, r\geq 1$
nous obtenons sa fonction g\'en\'eratrice 
\begin{eqnarray*}
&&\sum_\mu \frac{t^{|\mu|}}{z_\mu}\prod_{i\geq
1}\left[(1-(1-q)^i)\left(z+\sum_{j\geq
1}u^j\frac{(i)_j}{j!}p_j(Y)\right)\right]^{m_i(\mu)}\\
  & = & \prod_{i\geq 1}\sum_{m_i\geq =
0}\frac{t^{i\,m_i}}{m_i!\,i^{m_i}}
\left[(1-(1-q)^i)(z+\sum_{j\geq 1}((1-y_ju)^{-i}-1))\right]^{m_i}.
\end{eqnarray*}
Mais le dernier terme peut s'\'ecrire
\begin{eqnarray*}
&&\prod_{i\geq 1}\exp\left\{\left({t^i\over i}-{t^i\over 
i}(1-q)^i\right)
\left(z+\sum_{j\geq 1}((1-y_ju)^{-i}-1)\right)\right\}\\
&=&\left(1+{tq\over 1-t}\right)^z\prod_{j\geq 1}\exp\sum_{i\geq 1}
\left({t^i\over i}-{t^i\over i}(1-q)^i\right)\left({1\over 
(1-y_ju)^i}-1\right)\\
&=&\left(1+{tq\over 1-t}\right)^z
\prod_{j\geq 1}\left(1+{tq\over 1-t-y_ju}\right)\left(1+{tq\over 
1-t}\right)^{-1}.
\end{eqnarray*}
On pourrait trouver
des calculs similaires aux pr\'ec\'edents dans \cite{Las1} 
ou dans \cite{LL} en termes de $\lambda$-anneaux.


\end{document}